\documentclass[11pt]{article}%
\usepackage{amssymb}
\usepackage{amsfonts}
\usepackage{amsmath}
\usepackage{graphicx}%
\providecommand{\U}[1]{\protect\rule{.1in}{.1in}}
\begin{document}

\title{Estimating the covariance structure of heterogeneous SIS epidemics on networks }
\author{E. Cator\thanks{Radboud University Nijmegen, Faculty of Science, P.O. Box
9010, 6500 GL Nijmegen, The Netherlands; \emph{email}: E.Cator@science.ru.nl
}, \ H. Don\thanks{Amsterdam Free University, Faculty of Sciences, De
Boelelaan 1081, 1081 HV Amsterdam, The Netherlands, \emph{email}: h.don@vu.nl}
\ and P. Van Mieghem\thanks{Delft University of Technology, Faculty of EECS,
P.O. Box 5031, 2600 GA Delft, The Netherlands; \emph{email}:
P.F.A.VanMieghem@tudelft.nl }}
\date{21 July 2017}
\maketitle

\begin{abstract}
Heterogeneous Markovian Susceptible-Infected-Susceptible (SIS) epidemics with
a general infection rate matrix $\widetilde{A}$ are considered. Using a non-negative matrix factorization to approximate $\widetilde{A}$, we are able
to identify when a metastable state can be expected, and that the metastable
distribution, under certain conditions, will feature a normal distribution
with known expectation and covariance. Furthermore, we model a heterogeneous
Markovian SIS epidemic, that starts with a fraction of initially infected
nodes different from that in the metastable state, by approximating its
behaviour by a standard linear stochastic differential equation (SDE) in
sufficiently high dimensions. By exploiting the knowledge of the covariance
matrix from the SDE, we demonstrate significant accuracy improvements over the
first-order mean-field approximation NIMFA.

\end{abstract}

\section{Introduction}

Epidemic processes on a network can model an amazingly large variety of
real-world processes \cite{PVM_RMP_epidemics2014}, such as the spread of a
disease, a digital virus, a message in an on-line social network, an emotion,
the propagation of a failure or an innovation and other diffusion phenomena on
networks. We confine ourselves here to a particularly simple epidemic model,
the Markovian SIS epidemic process, that, as earlier
\cite{PVM_SIS_universality_tanh} argued, \textquotedblleft allows for the
highest degree of analytic treatment, which is a major motivation for the
continued effort towards its satisfactory understanding\textquotedblright.
Here and in addition to \cite{PVM_SIS_universality_tanh}, we justify that
quote. Whereas \cite{PVM_SIS_universality_tanh} was based on spectral
decomposition, here non-negative matrix factorization and clustering leads to
yet another powerful new analytic approximation, entirely different from
mean-field theory, but with an accuracy comparable to the well-established
mean-field models for large networks. Before outlining the idea of this paper,
we introduce definitions and notations.

We consider an unweighted, undirected graph $G$ containing a set $\mathcal{N}$
of $n$ nodes (also called vertices) and a set $\mathcal{L}$ of $L$ links (or
edges). The topology of the graph $G$ is represented by a symmetric $n\times
n$ adjacency matrix $A$ with elements $a_{ij}$. In an SIS epidemic process
\cite{Bailey_book,Anderson_May,Daley,Diekmann_Heesterbeek_Britton_boek2012,PVM_RMP_epidemics2014,Kiss_Miller_Simon2016}
on the graph $G$, the viral state of a node $i$ at time $t$ is specified by a
Bernoulli random variable $X_{i}\left(  t\right)  \in\{0,1\}$: $X_{i}\left(
t\right)  =0$, when node $i$ is healthy, but susceptible and $X_{i}\left(
t\right)  =1$, when node $i$ is infected at time $t$. A node $i$ at time $t$
can only be in one of these two states: \emph{infected}, with probability
$\Pr[X_{i}(t)=1]$ or \emph{healthy}, but susceptible to the infection, with
probability $\Pr[X_{i}(t)=0]=1-\Pr[X_{i}(t)=1]$. We assume that the curing
process for node $i$ is a Poisson process with rate $\delta_{i}$ and that the
infection rate over the directed link $\left(  i,j\right)  $ is a Poisson
process with rate $\beta_{ij}$. Only when node $i$ is infected, it can infect
each node $k$ of its healthy direct neighbors with rate $\beta_{ik}$. All
Poisson curing and infection processes are independent. This description
defines the continuous-time, Markovian \emph{heterogeneous} SIS epidemic
process on a graph $G$. We do not consider non-Markovian epidemics
\cite{PVM_nonMarkovianSIS_2013,PVM_nonMarkovianSIS_NIMFA_2013} and assume that
the infection characteristics in the graph, i.e. all curing and infection
rates, are independent of time.

Exploiting the important property $E\left[  X_{i}\right]  =\Pr\left[
X_{i}=1\right]  $ of a Bernoulli distribution, which enables to avoid
computations with the probability operator in favor of the easier, linear
expectation operator, the exact Markovian \emph{heterogeneous} SIS governing
equation \cite{PVM_secondorder_SISmeanfield_PRE2012,PVM_PAComplexNetsCUP} for
the infection probability of node $i$ is%
\begin{equation}
\frac{dE\left[  X_{i}\left(  t\right)  \right]  }{dt}=E\left[  -\delta
_{i}X_{i}\left(  t\right)  +\left(  1-X_{i}\left(  t\right)  \right)
\sum_{k=1}^{n}\beta_{ki}a_{ki}X_{k}\left(  t\right)  \right]
\label{governing_eq_heterogeneous_SIS}%
\end{equation}
When node $i$ is infected at time $t$, $X_{i}\left(  t\right)  =1$ and only
the first term on the right-hand side between the brackets $\left[  .\right]
$ affects the change in infection
probability with time $\frac{d\Pr\left[  X_{i}\left(  t\right)  =1\right]
}{dt}$ (left-hand side in (\ref{governing_eq_heterogeneous_SIS})). When node
$i$ is healthy, $X_{i}\left(  t\right)  =0$ and only the second term on the right-hand side between
the brackets $\left[  .\right]  $ increases the change in infection
probability with time by a rate $\sum_{k=1}^{n}\beta_{ki}a_{ki}X_{k}\left(
t\right)  $ due to all its infected, direct neighbors. We define the nodal
curing vector $\widetilde{\delta}=\left(  \delta_{1},\delta_{2},\ldots
,\delta_{n}\right)  $ and the weighted, infection rate adjacency matrix
$\widetilde{A}$ with element $\widetilde{a}_{ij}=\beta_{ij}a_{ij}$, which is
not necessarily symmetric. The governing equation
(\ref{governing_eq_heterogeneous_SIS}) can be evolved into a
\emph{heterogeneous} SIS Markov chain containing $2^{n}$ states
\cite{PVM_ToN_VirusSpread,PVM_EpsilonSIS_PRE2012,PVM_PAComplexNetsCUP}, and
although the Markov equations are linear and determine the joint probability
$\Pr\left[  X_{j_{1}}\left(  t\right)  =1,X_{j_{2}}\left(  t\right)
=1,\ldots,X_{j_{k}}\left(  t\right)  =1\right]  $ of each $k$-tuple $\left(
j_{1},j_{2},\ldots,j_{k}\right)  $ of nodal indices for $1\leq k\leq n$, the
solution of a set of $2^{n}$ linear equations is unfeasible for network sizes
$n>20$. For special types of graphs that feature symmetry, such as the
complete graph and the star \cite{PVM_MSIS_star_PRE2012}, fortunately, that
huge set of $2^{n}$ linear equations can be considerably reduced to a linear
in $n$ number of linear equations, that determine the probability that
$k\in\left\{  0,1,\ldots,n\right\}  $ nodes are infected at time $t$.

Our new idea here starts from symmetry as in the complete graph with equal
infection rates, where nodes are indistinguishable in the metastable state.
Szemer\'{e}di's regularity lemma, as explained in \cite{Diestel}, roughly
tells us that any graph can be partitioned into a bounded number of equal
subgraphs and the links among those subgraphs are fairly uniformly
distributed, like in random graphs. Here, we will combine both concepts: the
indistinguishability of the nodes and the inherent, topological partition
structure of a graph. First, we approximate the infection rate adjacency
matrix $\widetilde{A}$ by invoking non-negative matrix factorization (NMF)
\cite{Hannah_LAA1983}\ as%
\[
\widetilde{A}_{n\times n}\simeq\left(  W^{T}\right)  _{n\times k}H_{k\times n}%
\]
where $W$ and $H$ are $k\times n$ matrices, where in general $k<<n$. Earlier
in \cite{PVM_SIS_communityNetworks2014}, the influence of a community
structure on the spread of epidemics was studied. By invoking the concept of
equitable partitions and the quotient matrix \cite[p. 23-25]{PVM_graphspectra}%
, Bonaccorsi \emph{et al.} \cite{Bonaccorsi_SIAM2015} have formalized the
observations in \cite{PVM_SIS_communityNetworks2014} to clustering in
networks, that exhibit a relatively easily recognizable partitioning (such as
the interconnection of cliques). Based on graph partitioning and the isoperimetric
inequality, the universal mean-field framework (UMFF) \cite{PVM_UMFF} further
generalizes graph partitioning for mean-field epidemics and incorporates a
wide variety of different mean-field approximations such as
NIMFA\ \cite{PVM_N_intertwined_Computing2011} and the Heterogeneous Mean-Field
(HMF) \cite{PastorSatorras_Vespignani_PRE2001}. We believe that NMF is
generally applicable to networks, that do not immediately unveil an obvious
partition structure. The second step in our approach is to study the SIS
behaviour whenever the infection rate matrix $\widetilde{A}$ is of the form
$W^{T}H$. It is well known \cite{Lee_Nature1999} that NMF can be used to
perform a clustering of the $n$ nodes. This clustering or partitioning idea
transforms the SIS epidemics into a solvable stochastic diffusion process,
which is amenable not only to the analysis of the quasi-stationary or
metastable configuration, but even to the analysis of the time-evolution of
the SIS epidemics as shown in Section \ref{sec_continuous_process}.

Our method, called \textquotedblleft clustered SIS\textquotedblright, is
evaluated in Section \ref{sec_evaluation} on both a synthetic power law graph
as well as on a real-world network, derived from an airport network. Aimed to
encompass most mean-field approximations, UMFF \cite{PVM_UMFF} only provides
general bounds on the accuracy of mean-field approximations derived from the
isoperimetric inequality. UMFF does not propose mean-field accuracy
improvements techniques, in contrast to this paper. In particular, the
estimate of the covariance matrix $\Sigma_{\infty}$ in the metastable state,
deduced from the stochastic vector differential equation
(\ref{SDE_non_linear_D}) by linearization (Appendix
\ref{sec_exact_solution_SVDE}) is given in a closed form expression
(\ref{Sigma_inf_explicit_Sigma}) and its incorporation significantly improves
the accuracy of the first-order mean-field approximation NIMFA as demonstrated
in Section \ref{sec_application_real_world_networks} on a real-world network.
While earlier no specific information about the correlation structure
$E\left[  X_{i}X_{k}\right]  $ between any pair $\left(  i,k\right)  $ of
nodes was available (apart from the inequality $E\left[  X_{i}X_{k}\right]
\geq E\left[  X_{i}\right]  E\left[  X_{k}\right]  $ proved in
\cite{PVM_postive_SIScorrelations}), the common mean-field argument, assuming
independence and leading to the approximation $E\left[  X_{i}X_{k}\right]
\approx E\left[  X_{i}\right]  E\left[  X_{k}\right]  $, was a best possible
guess. Our estimate $\Sigma_{\infty}$ of the correlation structure thus adds
crucial new information to the theory of SIS epidemics on networks.

\section{Non-negative matrix factorization}

\label{sec_non_negative_matrix_factorization}In recent years, a lot of
attention has been devoted to NMF. Several good algorithms have been developed
to find the matrices $W$ and $H$ in ${\mathbb{R}}_{+}^{k\times n}$ that
minimise given objective functions such as the Frobenius norm $\Vert
\widetilde{A}-W^{T}H\Vert$ (see e.g. \cite[p. 234]{PVM_graphspectra}). A
particularity of our specific problem is that the diagonal of $\widetilde{A}$
is zero and irrelevant (there is no self-infection), which may lead to
different objective functions and possibly adapted algorithms for finding $W$
and $H$. In the case where $\widetilde{A}$ is symmetric, it is common practice
\cite{Vandaele_IEEE_TSP2016} to choose $H=W$, but we will not impose this
equality in general. In the sequel, we will assume that NMF algorithms will
supply us with the matrices $W$ and $H$ such that $W^TH$ approximates $\widetilde{A}$, 
except possibly on the diagonal.

After this first approximation, we continue by assuming that $\widetilde{A}=W^{T}H$ outside of the diagonal. Denote the $k\times n$ matrices
$W=\left[
\begin{array}
[c]{cccc}%
W_{1} & W_{2} & \cdots & W_{n}%
\end{array}
\right]  $ and $H=\left[
\begin{array}
[c]{cccc}%
H_{1} & H_{2} & \cdots & H_{n}%
\end{array}
\right]  $, where $W_{i}$ and $H_{i}$ are real $k\times1$ vectors. The
infection rate at which node $i$ infects node $j$ is then given by the inner
product of the vectors $W_{i}$ and $H_{j}$:
\begin{equation}
\widetilde{a}_{ij}=\beta_{ij}a_{ij}=W_{i}^{T}H_{j} \label{infection_rate_ij}%
\end{equation}
We therefore call $W_{i}$ the \textquotedblleft
infectiousness\textquotedblright\ of node $i$, and $H_{j}$ the
\textquotedblleft susceptibility\textquotedblright\ of node $j$. The final
characteristic of a node $i$ is its healing rate $\delta_{i}$. Hence, each
node $i$ is characterized by a vector $Z_{i}\in{\mathbb{R}}_{+}^{2k+1}$, with
\begin{equation}
Z_{i}=%
\begin{pmatrix}
\sqrt{n}\,W_{i}\\
\sqrt{n}\,H_{i}\\
\delta_{i}%
\end{pmatrix}
. \label{def_Zi}%
\end{equation}
The scaling by $\sqrt{n}$ of the two vectors $W_{i}$
and $H_{i}$ is useful later on in Section \ref{sec_continuous_process}. If two
nodes $i$ and $j$ have the same \textquotedblleft Z-vector\textquotedblright%
\ (so $Z_{i}=Z_{j}$), then these two nodes become indistinguishable in the SIS
process. Analogous to the complete graph, this means that the reduced process
where we only remember how many of those nodes with equal \textquotedblleft
Z-vector\textquotedblright\ are infected, is still a Markov process. If
$n>>k$, then there are many Z-vectors in ${\mathbb{R}}_{+}^{2k+1}$, implying
that many of the $n$ nodes are almost indistinguishable and such
indistinguishable nodes can be clustered and represented by a single $Z$-vector.

Let us choose $r$ clusters $C_{j}\subset\{1,\ldots,n\}$ of nodes such that (a)
their union $\cup_{j=1}^{r}C_{j}=\{1,\ldots,n\}$ comprises all the nodes in
the graph $G$ and (b) clusters are not overlapping, i.e. $C_{j}\cap
C_{l}=\emptyset$ if $j\neq l$, meaning that two different clusters do not
share a node. We define by $n_{j}$ the number of nodes in cluster $C_{j}$ and
we call $Y_{j}$ the \emph{cluster center,} defined as
\begin{equation}
Y_{j}=\frac{1}{n_{j}}\sum_{i\in C_{j}}Z_{i} \label{def_Yj}%
\end{equation}
which equals the average of all $Z$-vectors in cluster $j$. Ideally, we want
to choose the $r$ clusters in such a way that almost all nodes have their
$Z$-vector reasonably close to their corresponding cluster center $Y_{j}$.
Intuitively, the number $r$ of clusters must be sufficiently large such that
each cluster center $Y_{j}$ is a reasonable approximation for most nodes in
the cluster $C_{j}$. However, $r$ should not be too large, otherwise the
number of nodes in each cluster may not be large enough. It will turn out that 
this last requirement is not as essential as being close to the cluster center.\newline

Once the set of clusters $\left\{  C_{j}\right\}  _{1\leq j\leq r}$ is
determined, we propose a new \textquotedblleft clustered SIS\textquotedblright%
\ model, that approximates the original Markovian \emph{heterogeneous} SIS
process, by replacing the $Z_{i}$-vector of each node $i\in C_{j}$ by its
corresponding cluster center $Y_{j}$. In other words, we replace $Z_{i}$ by
\[
\tilde{Z}_{i}=\sum_{j=1}^{r}1_{\{i\in C_{j}\}}Y_{j},
\]
and any node $i$ belonging to the same cluster $C_{j}$ thus possesses the same
$\tilde{Z}_{i}=Y_{j}$. In our new approximate \textquotedblleft clustered
SIS\textquotedblright\ model, all nodes belonging to the same cluster are
indistinguishable. Hence, the random vector $N=(N_{1},\ldots,N_{r})$
containing the number $N_{j}$ of infected nodes (with $0\leq N_{j}\leq n_{j}$)
in each cluster $C_{j}$ is a continuous-time Markov process! We will study
this Markov process in more detail.

\section{The Markov process $N=(N_{1},\ldots,N_{r})$}

\label{sec_Markov_process}We introduce the following notation: we split up any
vector $y\in{\mathbb{R}}^{2k+1}$ according to the block vector structure of
$Z_{i}$,
\[
y=%
\begin{pmatrix}
y_{w}\\
y_{h}\\
y_{\delta}%
\end{pmatrix}
,
\]
where the vector $y_{w}\in{\mathbb{R}}^{k}$ corresponds to $\sqrt{n}$ times
the infectiousness, the vector $y_{h}\in{\mathbb{R}}^{k}$ corresponds to
$\sqrt{n}$ times the susceptibility and the scalar $y_{\delta}$ corresponds to
the healing rate. In particular, the cluster center is denoted by $Y_{j}=$ $%
\begin{pmatrix}
Y_{w,j}\\
Y_{h,j}\\
Y_{\delta,j}%
\end{pmatrix}
$.

The generator of the process $N\left(  t\right)  $ is determined by the
possible transition rates:
\begin{equation}
N_{j}\longrightarrow N_{j}+1\mbox{ at rate }\frac{1}{n}\sum_{l=1}^{r}%
(n_{j}-N_{j})N_{l}Y_{w,l}^{T}Y_{h,j} \label{eq:rateN1}%
\end{equation}
while%
\begin{equation}
N_{j}\longrightarrow N_{j}-1\mbox{ at rate }N_{j}Y_{\delta,j}.
\label{eq:rateN2}%
\end{equation}
In our new approximate \textquotedblleft clustered SIS\textquotedblright%
\ model, the number $N_{j}$ of infected nodes in cluster $C_{j}$ increases as
in (\ref{eq:rateN1}), because each of the $N_{l}$ infected nodes in cluster
$C_{l}$ can infect a healthy node of the total $n_{j}-N_{j}$ healthy nodes in
cluster $C_{j}$ at rate equal to $\frac{1}{n}Y_{w,l}^{T}Y_{h,j}$, which even
holds true for $l=j$.
The rate equation (\ref{eq:rateN2}) means that each of the $N_{j}$ infected
nodes in cluster $C_{j}$ can heal with rate $Y_{\delta,j}$. The number $N_{j}$
of infected nodes in cluster $C_{j}$ follows a birth-death process
\cite{PVM_PAComplexNetsCUP}, in which the birth rates depend upon the number
$N_{l}$ of infected nodes in any cluster $C_{l}$, including $C_{j}$ itself. In
fact, the new approximate \textquotedblleft clustered SIS\textquotedblright%
\ model for the vector $N$, that approximates the number of infected nodes in
$G$, consists of $r$ birth-death processes, coupled to each other through the
infection rates $\frac{1}{n}Y_{w,l}^{T}Y_{h,j}$ from cluster $C_{l}$ to
$C_{j}$.

Clearly, the state space of the continuous-time process $N\left(  t\right)
=(N_{1}(t),\ldots,N_{r}\left(  t\right)  )$, in which the integer random
variable $N_{j}\left(  t\right)  $ changes with time $t$, equals the product
$\bigotimes_{j=1}^{r}\left\{  0,1,2,\ldots,n_{j}\right\}  $ of possible
infection states per cluster $C_{j}$. We now embed this state space in
${\mathbb{R}}_{+}^{r}$ to allow the computation with real, rather than integer
$r\times1$ vectors. In particular, the metastable state of the new approximate
\textquotedblleft clustered SIS\textquotedblright\ model is obtained for a
real vector $N^{\infty}\in{\mathbb{R}}_{+}^{r}$ when, for each component
$N_{j}$, the rate of increase by one equals the rate of decrease by one. These
equalities are commonly called the balance equations for the steady state in a
continuous-time Markov process \cite[p. 213-214]{PVM_PAComplexNetsCUP}. The
vector $N^{\infty}$ represents the expected number of infected nodes per
cluster in meta-stable equilibrium. Equating the rate equations
(\ref{eq:rateN1}) and (\ref{eq:rateN2}) yields
\begin{equation}
\forall\ j\in\{1,\ldots,r\}\ :\ \frac{1}{n}\sum_{l=1}^{r}(n_{j}-N_{j}^{\infty
})N_{l}^{\infty}Y_{w,l}^{T}Y_{h,j}=N_{j}^{\infty}Y_{\delta,j}.
\label{eq:defNinfty}%
\end{equation}
which is an $r$-dimensional vector equation, whose solution is shown to exist
below. If $k<r$, we can solve these equations in the following way: define the
vector $V\in{\mathbb{R}}^{k}$ by
\begin{equation}
V=\frac{1}{n}\sum_{l=1}^{r}N_{l}^{\infty}Y_{w,l}. \label{def_V}%
\end{equation}
Then \eqref{eq:defNinfty} becomes
\[
\forall\ j\in\{1,\ldots,r\}\ :\ (n_{j}-N_{j}^{\infty})V^{T}Y_{h,j}%
=N_{j}^{\infty}Y_{\delta,j}.
\]
from which it follows that
\begin{equation}
N_{j}^{\infty}=n_{j}\cdot\frac{V^{T}Y_{h,j}}{V^{T}Y_{h,j}+Y_{\delta,j}}.
\label{eq:Ninfty}%
\end{equation}
It remains to determine the vector $V$. Substituting (\ref{eq:Ninfty}) into
the definition (\ref{def_V}) of $V$ yields the recursion
\begin{equation}
V=\frac{1}{n}\sum_{l=1}^{r}n_{l}\cdot\frac{Y_{w,l}Y_{h,l}^{T}V}{Y_{h,l}%
^{T}V+Y_{\delta,l}}. \label{eq:V}%
\end{equation}
So we have reduced the $r$-dimensional equation (\ref{eq:defNinfty}) to a
$k$-dimensional equation. This $k$-dimensional vector equation (\ref{eq:V}) in
$V$ can numerically be solved efficiently, even for relatively large $k$, provided a
solution exists.

\subsection{Existence of the metastable state}

\label{sec_metastable_state}We will now write the equation
\eqref{eq:defNinfty} for the metastable state in matrix form. Let us first
define $\mathbf{n}$ as the vector $(n_{1},\ldots,n_{r})$ with the number of
nodes in each cluster and denote the $k\times r$ matrices
\[
Y_{w}=\left[
\begin{array}
[c]{ccc}%
Y_{w,1} & \cdots & Y_{w,r}%
\end{array}
\right]  \text{ and }Y_{h}=\left[
\begin{array}
[c]{ccc}%
Y_{h,1} & \cdots & Y_{h,r}%
\end{array}
\right]  .
\]
In addition, we use the following notation for a vector $a\in{\mathbb{R}}^{p}%
$:
\[
\mathrm{diag}(a)=%
\begin{pmatrix}
a_{1} &  & \\
& \ddots & \\
&  & a_{p}%
\end{pmatrix}
.
\]
Then, \eqref{eq:defNinfty} becomes%
\[
\frac{1}{n}\mathrm{diag}(\mathbf{n})Y_{h}^{T}Y_{w}N^{\infty}-\frac{1}%
{n}\mathrm{diag}(N^{\infty})Y_{h}^{T}Y_{w}N^{\infty}=\mathrm{diag}(Y_{\delta
})N^{\infty}%
\]
and%
\begin{align*}
0  &  =\left\{  \frac{1}{n}\mathrm{diag}(\mathbf{n})Y_{h}^{T}Y_{w}%
-\mathrm{diag}(Y_{\delta})\right\}  N^{\infty}-\frac{1}{n}\mathrm{diag}%
(N^{\infty})Y_{h}^{T}Y_{w}N^{\infty}\\
&  =\mathrm{diag}(Y_{\delta})\left\{  \frac{1}{n}\mathrm{diag}(\mathbf{n}%
/Y_{\delta})Y_{h}^{T}Y_{w}-I\right\}  N^{\infty}-\frac{1}{n}\mathrm{diag}%
(N^{\infty})Y_{h}^{T}Y_{w}N^{\infty}%
\end{align*}
After denoting the non-negative $r\times r$ matrix
\[
\overline{A}=\frac{1}{n}\mathrm{diag}(\mathbf{n}/Y_{\delta})Y_{h}^{T}Y_{w},
\]
the balance equation \eqref{eq:defNinfty} finally becomes%
\begin{equation}
0=\mathrm{diag}(Y_{\delta})(\overline{A}-I)N^{\infty}-\frac{1}{n}%
\mathrm{diag}(N^{\infty})Y_{h}^{T}Y_{w}N^{\infty}
\label{matrix_eq_metastable_state}%
\end{equation}
As in \cite{PVM_heterogeneous_virusspread}, we can prove that the metastable
state matrix equation (\ref{matrix_eq_metastable_state}) has a strictly
positive solution in $N^{\infty}$ if the $r\times r$-matrix $\overline{A}$ has
a positive eigenvalue greater than 1. Indeed, let us assume that all healing
rates $\delta_{i}=\delta$ are the same, then
\[
\overline{A}=\frac{1}{\delta}\frac{1}{n}\mathrm{diag}(\mathbf{n})Y_{h}%
^{T}Y_{w}.
\]
Furthermore, suppose that each $Z_{i}$ is very close to its cluster center,
for example, if $r$ is close to $n$. Next, define the map $\mathrm{Cl}%
:\{1,\ldots,n\}\rightarrow\{1,\ldots,r\}$, such that $\mathrm{Cl}(i)=j$
precisely when $i\in C_{j}$. Since the $n\times n$ infection rate adjacency
matrix $\widetilde{A}^{T}=H^{T}W$, we can write its elements using the map
$\mathrm{Cl}$ as
\[
\widetilde{a}_{ii^{\prime}}^{T}=H_{i}^{T}W_{i^{\prime}}=\frac{1}{n}(Y_{h}%
^{T}Y_{w})_{\mathrm{Cl}(i)\mathrm{Cl}(i^{\prime})}.
\]
Hence, if $\overline{v}\in{\mathbb{R}}^{r}$ is the positive eigenvector of the
non-negative matrix $\overline{A}$ corresponding to the eigenvalue $\lambda
>1$, then, after defining $v\in{\mathbb{R}}^{n}$ by
\[
v_{i}=\overline{v}_{\mathrm{Cl}(i)},
\]
we arrive at
\[
\frac{1}{\delta}(\widetilde{A}^{T}v)_{i}=\frac{1}{\delta}\sum_{i^{\prime}%
=1}^{n}\widetilde{a}_{ii^{\prime}}^{T}v_{i^{\prime}}=\frac{1}{\delta}%
\sum_{j=1}^{r}\frac{n_{j}}{n}(Y_{h}^{T}Y_{w})_{\mathrm{Cl}(i)j}\overline
{v}_{j}=(\overline{A}\overline{v})_{\mathrm{Cl}(i)}=\lambda v_{i},
\]
which demonstrates that $\lambda$ is also an eigenvalue of $\widetilde{A}%
^{T}/\delta$. This means that a metastable solution for $N^{\infty}$ exists if
the largest eigenvalue $\lambda$ of $\widetilde{A}^{T}/\delta$ is larger than
1, which is exactly the same condition that NIMFA
\cite{PVM_heterogeneous_virusspread} predicts in this case!

\section{Continuous approximation of the process $N$}

\label{sec_continuous_process}In this section, we show how our
\textquotedblleft clustered SIS\textquotedblright\ Markov process $N$, with
the transition rates given in \eqref{eq:rateN1} and \eqref{eq:rateN2}, can in
turn be approximated by a Markov process with a continuous state space in a
classical way, namely by approximating a difference of Poisson processes by a
Brownian motion with drift \cite{Harrison}.

We center and rescale the process $N$ and we define the deviations process $D$
by
\begin{equation}
D=\frac{N-N^{\infty}}{\sqrt{n}} \label{def_D}%
\end{equation}
where the $\sqrt{n}$-scaling will become clear later. The transition rates for
$D$ directly follow from \eqref{eq:rateN1} and \eqref{eq:rateN2},
\begin{align}
D_{j}  &  \longrightarrow D_{j}+1/\sqrt{n}\mbox{ at rate }\frac{1}{n}%
\sum_{l=1}^{r}(n_{j}-N_{j}^{\infty}-n^{1/2}D_{j})(n^{1/2}D_{l}+N_{l}^{\infty
})Y_{w,l}^{T}Y_{h,j}\label{eq:rateD1}\\
D_{j}  &  \longrightarrow D_{j}-1/\sqrt{n}\mbox{ at rate }(N_{j}^{\infty
}+n^{1/2}D_{j})Y_{\delta,j}. \label{eq:rateD2}%
\end{align}

We now suppose that the number of nodes $n$ in the graph $G$ is large. At some
time $t$, there will be many infection and healing events in a relatively
small time interval $[t,t+h]$. If $h$ is so small that the transition rates
can be assumed to be constant during the time interval $[t,t+h]$ (this means
that the changes in $D$ during that time interval are negligible), then the
distribution of the number $I_{j}(h)$ of nodes in cluster $C_{j}$ that become
infected during the time interval $[t,t+h]$ and the number $H_{j}(h)$ of nodes
in cluster $C_{j}$ that heal in that time interval can be determined. To a
good approximation, the random variables $I_{j}(h)$ and $H_{j}(h)$ for each
cluster $1\leq j\leq r$ will all be independent and Poisson distributed. We
can explicitly calculate the expectation of these random variables,
conditioned on $D(t)$:
\[
{\mathbb{E}}\left[  I_{j}(h)\mid D(t)\right]  =\frac{h}{n}\sum_{l=1}^{r}%
(n_{j}-N_{j}^{\infty}-n^{1/2}D_{j})(n^{1/2}D_{l}+N_{l}^{\infty})Y_{w,l}%
^{T}Y_{h,j}%
\]
and
\[
{\mathbb{E}}\left[  H_{j}(h)\mid D(t)\right]  =h(N_{j}^{\infty}+n^{1/2}%
D_{j})Y_{\delta,j}.
\]
Next, we calculate the expectation of the vector $dD(h)=D(t+h)-D(t)=n^{-1/2}%
(I(h)-H(h))$. The $j$-th component is
\begin{align*}
{\mathbb{E}}\left[  dD_{j}(h)\mid D(t)\right]   &  =n^{-\frac{1}{2}%
}(\mathbb{E}[I_{j}(h)\mid D(t)]-\mathbb{E}[H_{j}(h)\mid D(t)])\\
&  =hn^{-\frac{1}{2}}\left(  \frac{1}{n}\sum_{l=1}^{r}(n_{j}-N_{j}^{\infty
}-n^{1/2}D_{j})(n^{1/2}D_{l}+N_{l}^{\infty})Y_{w,l}^{T}Y_{h,j}-(N_{j}^{\infty
}+n^{1/2}D_{j})Y_{\delta,j}\right)  .
\end{align*}
Using the metastable state equation \eqref{eq:defNinfty} yields cancellation
of the main order terms:
\[
{\mathbb{E}}\left[  dD_{j}(h)\mid D(t)\right]  =h\left(  \frac{1}{n}\sum
_{l=1}^{r}\left(  (n_{j}-N_{j}^{\infty}-n^{1/2}D_{j})D_{l}-D_{j}N_{l}^{\infty
}\right)  Y_{w,l}^{T}Y_{h,j}-D_{j}Y_{\delta,j}\right)  .
\]
Conditionally on $D\left(  t\right)  $, the components of $dD(h)$ are
independent, so that we only need to determine the variances. Since the
variance of a Poisson distribution is equal to the mean, we add the
expectations of $I_{j}(h)$ and $H_{j}(h)$,
\begin{align*}
\mathrm{Var}\left[  dD_{j}(h)\mid D(t)\right]   &  =\mathrm{Var}\left[
n^{-\frac{1}{2}}(I_{j}(h)-H_{j}(h))\mid D(t)\right]  =\frac{1}{n}%
(\mathbb{E}[I_{j}(h)\mid D(t)]-\mathbb{E}[H_{j}(h)\mid D(t)])\\
&  =\frac{h}{n}\left(  \frac{1}{n}\sum_{l=1}^{r}(n_{j}-N_{j}^{\infty}%
-n^{1/2}D_{j})(n^{1/2}D_{l}+N_{l}^{\infty})Y_{w,l}^{T}Y_{h,j}+(N_{j}^{\infty
}+n^{1/2}D_{j})Y_{\delta,j}\right)  .
\end{align*}
and again introducing the metastable state equation \eqref{eq:defNinfty}, we
find%
\[
\mathrm{Var}\left[  dD_{j}(h)\mid D(t)\right]  =hn^{-\frac{1}{2}}\left(
\frac{1}{n}\sum_{l=1}^{r}\left(  (n_{j}-N_{j}^{\infty}-n^{1/2}D_{j}%
)D_{l}-D_{j}N_{l}^{\infty}\right)  Y_{w,l}^{T}Y_{h,j}+(2n^{-\frac{1}{2}}%
N_{j}^{\infty}+D_{j})Y_{\delta,j}\right)  .
\]
The addition of the expectations of $I_{j}(h)$ and $H_{j}(h)$ explains why
there is no cancellation of the main order terms of the variance, leading to
the term $2hN_{j}^{\infty}Y_{\delta,j}/n$. Compared to the expectation, there
is also an extra factor $n^{-\frac{1}{2}}$ in the variance. It is convenient
to write these equations for each $1\leq j\leq r$ in matrix notation as in
Section \ref{sec_metastable_state},
\[
{\mathbb{E}}\left[  dD(h)\mid D(t)\right]  =h\left(  \frac{1}{n}%
\mathrm{diag}(\mathbf{n}-N^{\infty}-n^{1/2}D)Y_{h}^{T}Y_{w}D-\frac{1}%
{n}\mathrm{diag}(Y_{h}^{T}Y_{w}N^{\infty})D-\mathrm{diag}(Y_{\delta})D\right)
.
\]
We introduce the $r\times r$-matrix $M$,
\begin{equation}
M=\frac{1}{n}\mathrm{diag}(\mathbf{n}-N^{\infty})Y_{h}^{T}Y_{w}-\frac{1}%
{n}\mathrm{diag}(Y_{h}^{T}Y_{w}N^{\infty}). \label{def_M}%
\end{equation}
and obtain
\begin{equation}
{\mathbb{E}}\left[  dD(h)\mid D(t)\right]  )=h\left(  (M-\mathrm{diag}%
(Y_{\delta}))D-\frac{1}{\sqrt{n}}\mathrm{diag}(D)Y_{h}^{T}Y_{w}D\right)  .
\label{eq:drift}%
\end{equation}
For the conditional variance of $dD(h)$, we find a similar formula:
\begin{equation}
\mathrm{Var}\left[  dD(h)\mid D(t)\right]  =h\left(  \frac{1}{\sqrt{n}%
}(M+\mathrm{diag}(Y_{\delta}))D-\frac{1}{n}\mathrm{diag}(D)Y_{h}^{T}%
Y_{w}D+\frac{2}{n}\mathrm{diag}(Y_{\delta})N^{\infty}\right)  .
\label{eq:variance}%
\end{equation}

For large $n$, it is natural to consider $D$ as a continuous process. We have
shown that $dD(h)$ is distributed as the scaled difference of independent
Poisson random variables, so that $dD(h)$ can be approximated by a normal
distribution. Equations \eqref{eq:drift} and \eqref{eq:variance} then suggest
to model $D$ as the solution of the following stochastic vector differential
equation (SVDE):
\begin{align}
dD(t)  &  =\left(  (M-\mathrm{diag}(Y_{\delta}))-\frac{1}{\sqrt{n}%
}\mathrm{diag}(D(t))Y_{h}^{T}Y_{w}\right)  D(t)dt+\nonumber\\
&  \ \ +\mathrm{diag}\left(  \sqrt{\frac{1}{\sqrt{n}}(M+\mathrm{diag}%
(Y_{\delta}))D-\frac{1}{n}\mathrm{diag}(D)Y_{h}^{T}Y_{w}D+\frac{2}%
{n}\mathrm{diag}(Y_{\delta})N^{\infty}}\,\right)  dB(t).
\label{SDE_non_linear_D}%
\end{align}
where $B(t)$ is an $r$-dimensional standard Brownian motion, and where we have
used the notation for vectors $a\in{\mathbb{R}}_{+}^{p}$,
\[
\sqrt{a}=(\sqrt{a_{1}},\sqrt{a_{2}},\ldots,\sqrt{a_{p}})
\]
Equation (\ref{SDE_non_linear_D}) can be used to efficiently simulate the
time-evolution of $D$. From such simulation, not only information about the
metastable distribution can be obtained, but also about the relaxation time,
for example. Since the SVDE in (\ref{SDE_non_linear_D}) is non-linear in $D$,
an exact analysis is hard. However, for large $n$, we can make additional
approximations as shown in Appendix \ref{sec_linearizing_SVDE}, that lead us
to conclude, provided the largest eigenvalue $\lambda$ of $A^{T}/\delta$ is
larger than $1$, that the metastable distribution of the vector of infected
nodes
\[
N=N^{\infty}+\sqrt{n}D
\]
equals an $r$-dimensional normal distribution, with expectation equal to
$N^{\infty}$ (since ${\mathbb{E}}\left[  D\right]  =0$), and covariance matrix
$n\Sigma_{\infty}$ specified in (\ref{Sigma_inf_explicit_Sigma}).

\section{Evaluation of the \textquotedblleft clustered SIS\textquotedblright%
\ model}

\label{sec_evaluation}

\subsection{Choosing the number of clusters $r$}

When applying our method to a particular example, we face two problems: we
have to choose both the dimension $k$ of our non-negative matrix factorization
of $\widetilde{A}$ and the number $r$ of clusters. The first problem is hard,
and it depends on how well $\widetilde{A}$ can be approximated by the product
of two low rank $k\times n$ matrices $W$ and $H$. Choosing the number $r$ of
clusters is less delicate: if we choose a lot of clusters, but a subset of the
cluster centers are very close together, then the clusters themselves become
indistinguishable. This means that the total number of infected nodes in the
union of the respective subset of clusters behaves as though it were one
cluster, also in the approximation of the SVDE. The argument guarantees a
certain stability: choosing more clusters does not really change the predicted
behavior of a fixed subset of nodes. We therefore recommend choosing a lot of
clusters (maybe even $r=n$), since the bias effect is less severe, whereas the
effect on the variance is limited, especially for larger collections of nodes.
Clearly, on a micro-level (when looking at a few nodes), the approximation
will not get much better.

\subsection{Application to a simulated network}

A random network with a power law degree distribution $\mathrm{Pr}(D_i>x)\sim
x^{-2}$, where $D_i$ is the degree of node $i$ in the network, is generated
by the configuration model. Only the largest connected component of the
generated network is used, consisting of $n=9994$ nodes and with a maximum
degree of 1229 and an average degree of 15. The infection rate matrix
$\widetilde{A}$ is symmetric, with infection rate $\beta_{ij}=1$ between all
pairs $\left(  i,j\right)  $ of connected nodes. We have made a one-dimensional
($k=1$) factorisation of $\widetilde{A}$, invoking a fast iterative algorithm
that minimises $\sum_{i=1}^{n}\sum_{j=1}^{n}\omega_{ij}(\widetilde{a}%
_{ij}-W_{i}^{T}H_{j})^{2}$. We take weights $\omega_{ij}$ of the form
$\omega_{ii}=0$ (the diagonal of $\widetilde{A}$ is irrelevant) and
$\omega_{ij}=e^{\lambda\widetilde{a}_{ij}}$ ($i\neq j$), where we can still
tune the parameter $\lambda\in{\mathbb{R}}$. In the current example, $\widetilde{A}$ is a 0-1-matrix, so for $i\neq j$ there are only two possible weights: $w_{ij}=1$ if $\tilde a_{ij}=0$ and $w_{ij}=e^\lambda$ if $\tilde a_{ij}=1$.  The parameter $\lambda$ allows us to control the
\textquotedblleft best\textquotedblright\ approximation $W^{T}H$ for
$\widetilde{A}$, since we have no theoretical results guiding us in this approximation. Here we choose $\lambda$ in such a way
that the NIMFA approximation of the expected total number of infected nodes for the network with matrix $W^TH$ is (approximately) equal to the original NIMFA approximation. 

\begin{figure}[ptb]%
\[
\includegraphics[bb=0 0 850 650, height = 0.4\textheight, width = 0.8\textwidth,clip=]{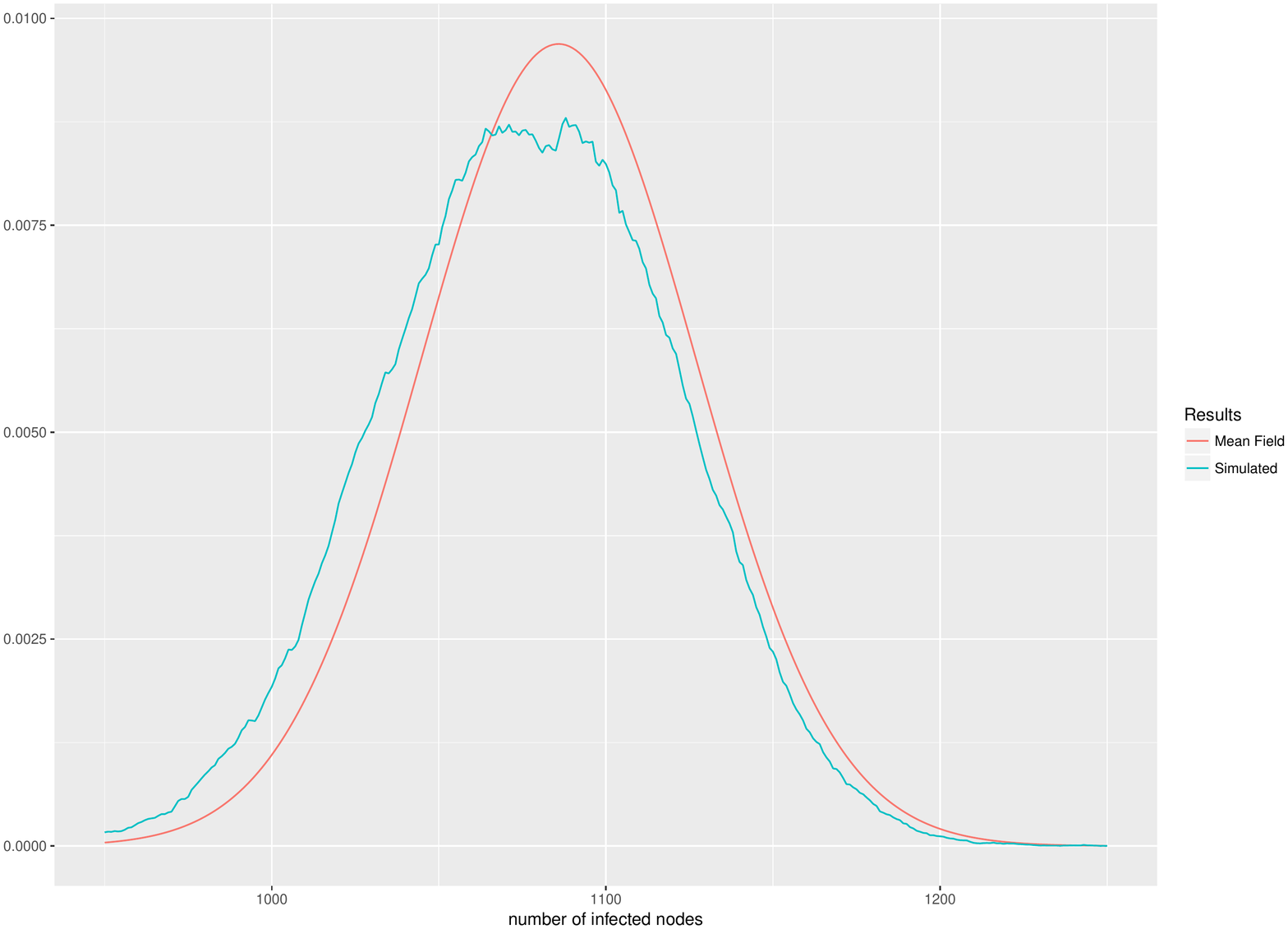}
\]
\caption{Distribution of number of infected nodes in the simulated network in
the metastable state.}%
\label{simul_cdf}%
\end{figure}

In the one dimensional approximation of the symmetric matrix $\widetilde{A}$, it comes out of the minimization procedure that also $W^TH$ will be symmetric, i.e.
we can choose the $n\times1$ matrix $W$ equal to $H$ and we base our
clustering solely on $W$. We choose $r=100$ clusters, where we put a few large
values of $W$ in separate clusters. The metastable state number of infected
nodes in our \textquotedblleft clustered SIS\textquotedblright\ model is
computed from (\ref{eq:Ninfty}) and (\ref{eq:V}), and the covariance matrix for all $100$ clusters is approximated by $n\Sigma_\infty$, specified in \eqref{Sigma_inf_explicit_Sigma}. Figure \ref{simul_cdf} shows
how well our 1-dimensional approximation is able to estimate the distribution
of the total number of infected nodes, which we obtained via simulation of the
entire network with infection matrix $\widetilde{A}$. The simulations started
with all nodes infected. Only when the metastable regime is approximately
reached, the simulation data is collected for a long time to accurately obtain
the metastable distribution of the infected nodes.
The accuracy of the $k=1$ approximation is also illustrated by the
following table. From the simulated distribution in the metastable state, the
expectation and standard deviation of the number of infected nodes in the
graph are computed in the first column and compared with the mean
$u^{T}N^{\infty}$ in (\ref{eq:Ninfty}), where $u$ is the all-one vector, and
the standard deviation $\sqrt{n}\cdot\sqrt{u^{T}\Sigma_{\infty}u}$ in
(\ref{Sigma_inf_explicit_Sigma}) in the last column. Note that this standard deviation depends on the clustering and is also used for plotting the curve in Figure \ref{simul_cdf}, as the mean field approximation only gives the center of the distribution. 
\[%
\begin{tabular}
[c]{|c|c|c|}\hline
& Simulated & Mean Field\\\hline
Expectation & 1076.0 & 1085.9\\\hline
Standard Deviation & 43.7 & 41.16\\\hline
\end{tabular}
\ \ \ \ \ \
\]
NIMFA \cite[Chapter 17]{PVM_PAComplexNetsCUP} provides an upper bound of viral
infection probability of a node. This causes the expectation to be overestimated by NIMFA. This overestimation is by construction still present in the curve in Figure \ref{simul_cdf}, but the estimation of the variance is rather good.

\begin{figure}[ptb]%
\[
\includegraphics[bb=0 0 850 650, height = 0.4\textheight, width = 0.8\textwidth,clip=]{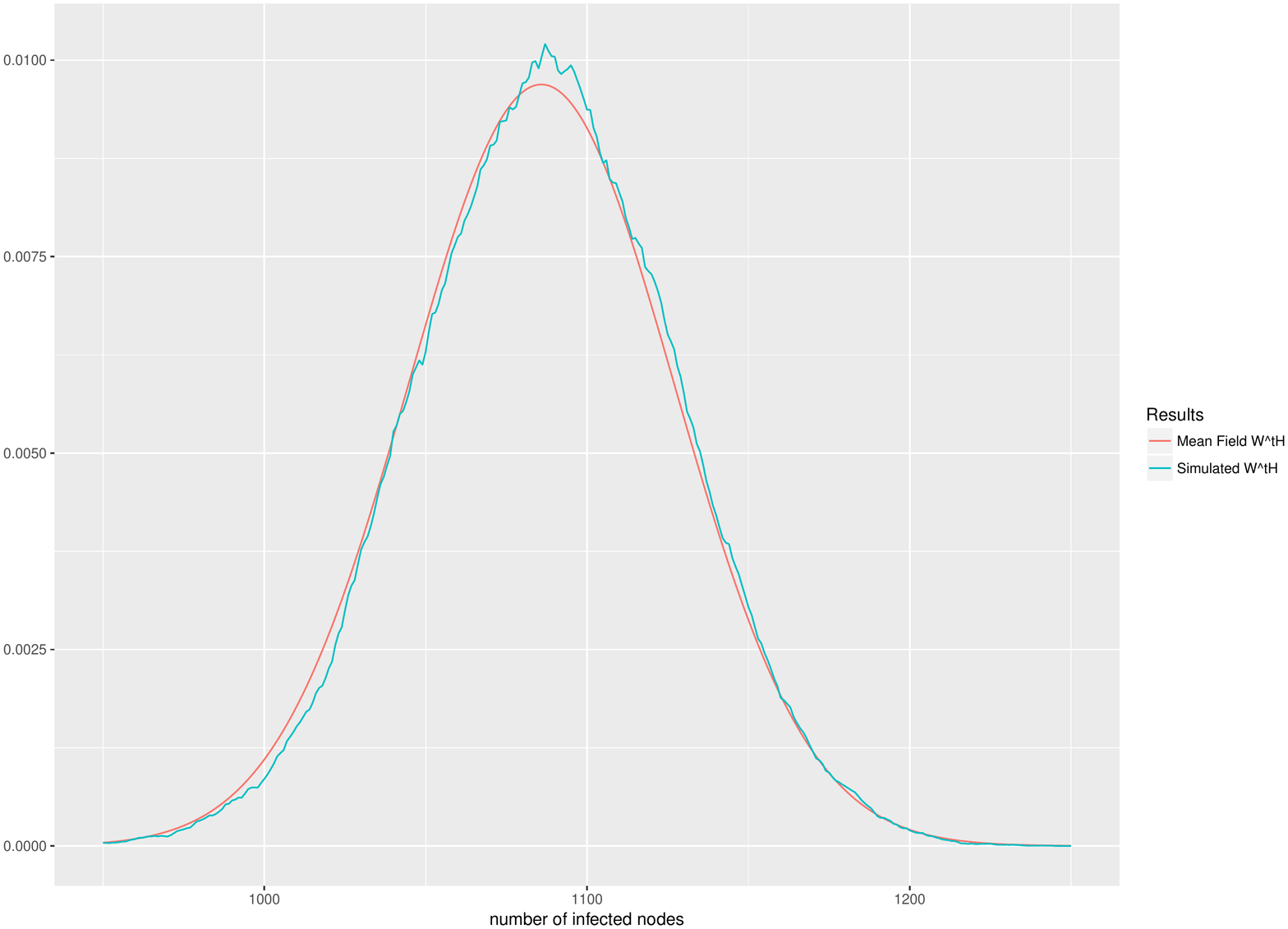}
\]
\caption{Distribution of number of infected nodes in the $W^{T}H$-network}%
\label{WtH_cdf}%
\end{figure}

Our approach makes several approximations: (a) the infection rate matrix
$\widetilde{A}$ is approximated by $W^{T}H$, which we further call the
$W^{T}H$-network, (b) clusters are chosen and (c) finally, the dynamics in
each cluster is described by a continuous Markov process $N$ as explained in
Section \ref{sec_Markov_process}. To assess the effects of these
approximations, we have also simulated directly from the $W^{T}H$-network in
the same way as before, but now using $W^{T}H$ as infection rate matrix, to
compare this with our predicted results. The results are astonishingly
accurate, as demonstrated in Figure \ref{WtH_cdf}. 
The predicted curve is the same as in Figure \ref{simul_cdf} and is very close to the $W^TH$-simulation. This means that for the $W^TH$-network the clustering and continuous approximation work really well; in Figure \ref{WtH_cdf} there is even no visible NIMFA bias. Notably, even the one-dimensional approach ($k=1$) is capable to approximate the variance well. 

The main issue that remains is the bias in the prediction for the original network. In the next paragraph, we will present a method to correct for this NIMFA bias, that, unfortunately, requires a high number of clusters and a high-dimensional approximation, which was unfeasible to compute for this large network.

\subsection{Application to a real-world network}

\label{sec_application_real_world_networks}We consider the airline network
(with $n=3425$ nodes and $L=15358$ links), where the nodes are the airports
and the infection rate $\beta_{ij}$ along a link is proportional to the number
of flights between the two airports. Qu and Wang \cite{Bo_Huijuan_TNSE2017}
have constructed this network and its infection rates from the dataset of
openflights\footnote{https://openflights.org/data.html}. Many of the
asymmetric infection rates ($\beta_{ij}\neq\beta_{ji}$) are zero; thus many
nodes\footnote{In fact, there are precisely 348 nodes with zero in- and out-degree.
Those nodes have been omitted by Qu and Wang \cite{Bo_Huijuan_TNSE2017}. Our
approach can handle such zero in- and zero out-degree nodes.} have zero
infectiousness (zero out-degree) or zero susceptibility (zero in-degree). This
means that the susceptibility and the infectiousness of a node are actually
different. 

\begin{figure}[ptb]%
\[
\includegraphics[bb=0 0 850 650, height = 0.4\textheight, width = 0.8\textwidth,clip=]{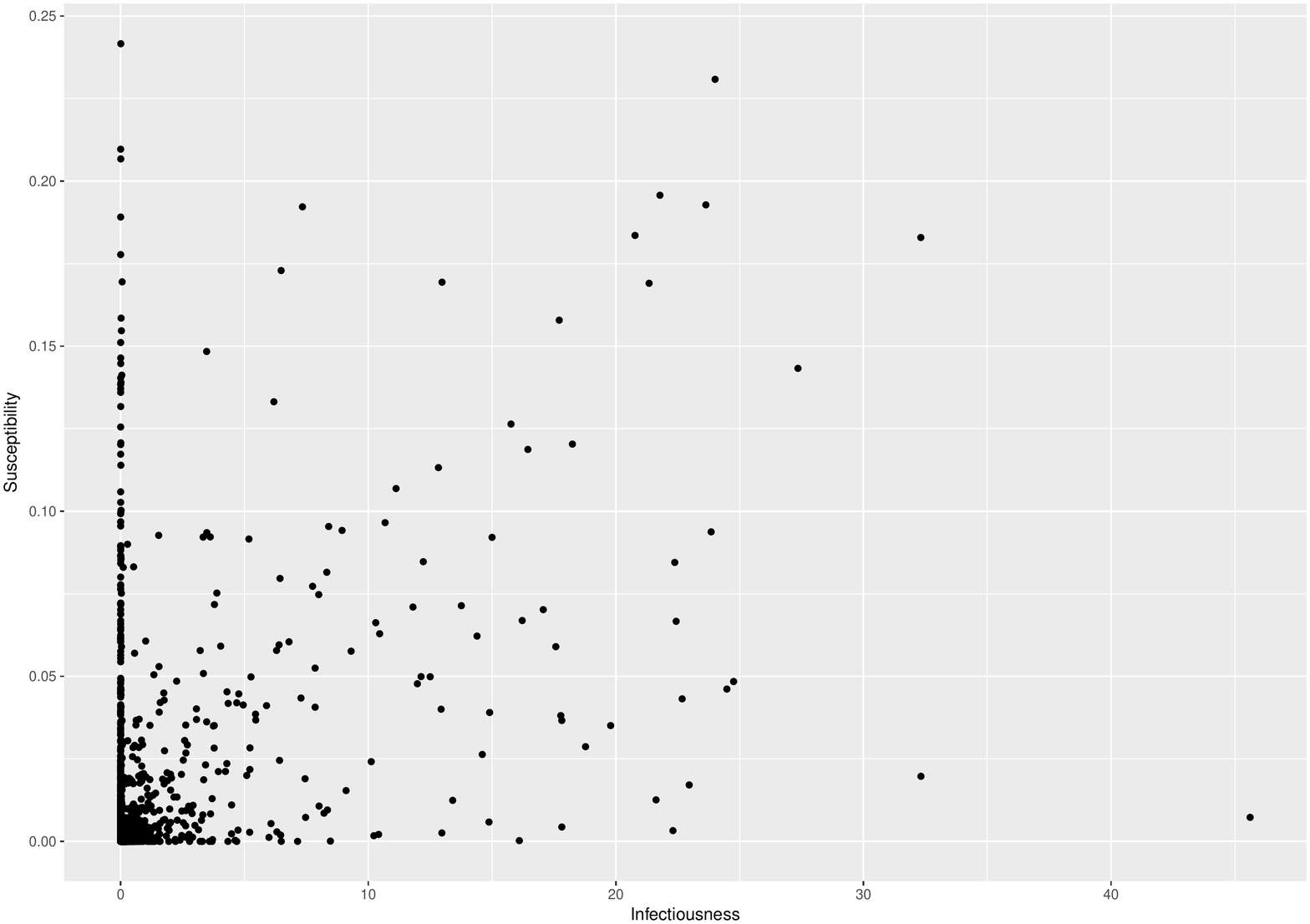}
\]
\caption{Scatter plot of $W$ and $H$ in the air weight network}%
\label{WHplot}%
\end{figure}

As before, we start by making a one dimensional factorization, resulting in a
vector $W$ and $H$, where we force the NIMFA approximations for the original and factorized matrix to be approximately equal. Figure \ref{WHplot} gives an
impression of $W$ and $H$. Clustering is not so obvious. Since the network is relatively
small, we decided to choose $n$ clusters: each node becomes a separate
cluster. However, the expected prevalence for each node (the vector
$N^{\infty}$) and the estimated covariance matrix $\Sigma_{\infty}$ are not
very impressive:
\[%
\begin{tabular}
[c]{|c|c|c|}\hline
& Simulated & Mean Field\\\hline
Expectation & 1105.8 & 1135.1\\\hline
Standard Deviation & 29.4 & 16.45\\\hline
\end{tabular}
\ \ \ \ \
\]
There is a sizeable overestimation of the expectation and a severe
underestimation of the variance. These differences are caused by the fact that in this case $\widetilde{A}$ is rather poorly approximated by the rank-one matrix $W^{T}H$. Apparently in some cases a one-dimensional factorization works pretty well, but for the current example clearly not. Additional research is needed to be able to predict performance of low-dimensional factorizations from network properties. 

\begin{figure}[ptb]%
\[
\includegraphics[bb=0 0 850 650, height = 0.4\textheight, width = 0.8\textwidth,clip=]{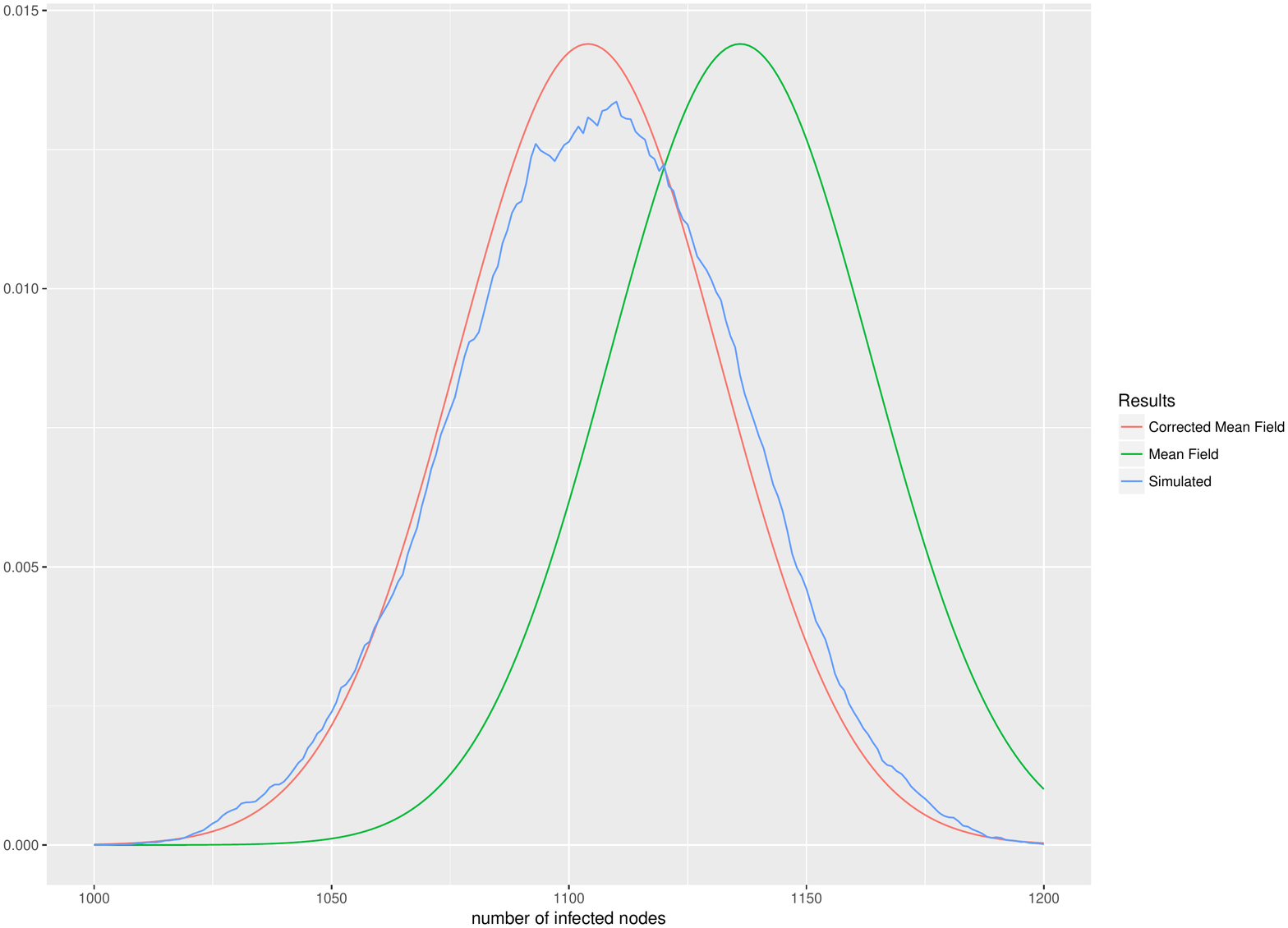}
\]
\caption{Distribution of number of infected nodes in Air Weight network}%
\label{AirWeight_cdf}%
\end{figure}

\begin{figure}[ptb]%
\[
\includegraphics[bb=0 0 850 650, height = 0.4\textheight, width = 0.8\textwidth,clip=]{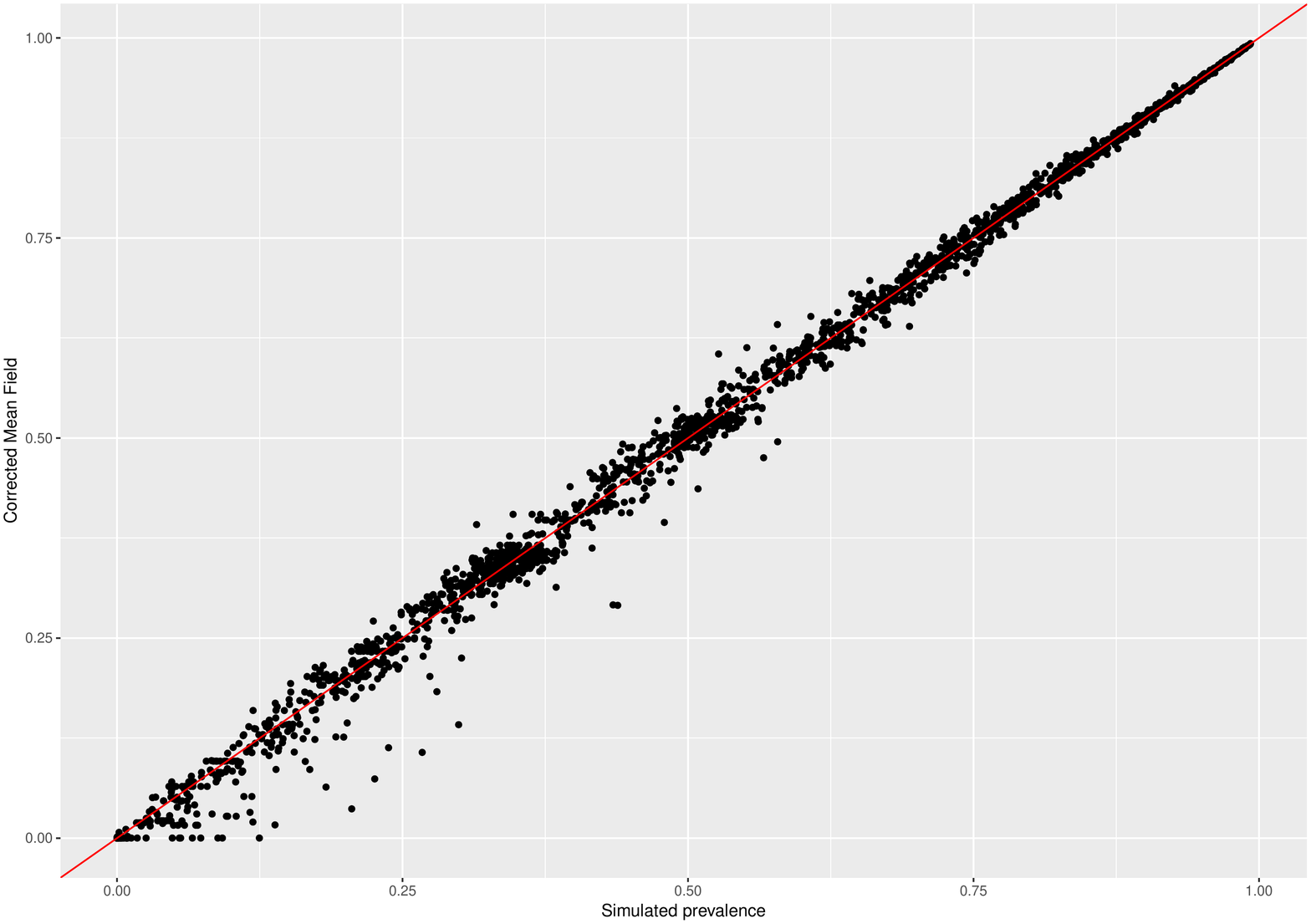}
\]
\caption{Infection probability of each node in Air Weight network}%
\label{AirWeight_exp}%
\end{figure}

To improve our results for this network, we need to increase the dimension $k$. In fact, since the network is small, we do not need to approximate $\widetilde{A}$ with a lower rank factorization at all: we simply use $\widetilde{A}$ itself in our
formula's (for example, choose $H=\widetilde{A}$ and $W=I$ the identity
matrix, and therefore $k=n$). Furthermore, each node in the network consists
of its own cluster. We also perform a correction on the NIMFA expectations
$N^{\infty}$ in the following way. NIMFA upper bounds the infection
probability per node: in governing equation
(\ref{governing_eq_heterogeneous_SIS}), NIMFA replaces $E\left[  X_{i}%
X_{k}\right]  $ by $E\left[  X_{i}\right]  E\left[  X_{k}\right]  $, ignoring
the positive covariance between $X_{i}$ and $X_{k}$. However, the matrix
$\Sigma_{\infty}$ provides an indication of the covariance between any two
nodes! We now use the following procedure: after determining $N^{\infty}$
using equations (\ref{eq:Ninfty}), we determine $\Sigma_{\infty}$ using
(\ref{Sigma_inf_explicit_Sigma}). Then, we determine a \emph{corrected vector}
$\widehat{N}^{\infty}$ of cluster expectations by adjusting the cluster rate
equations (\ref{eq:rateN1}) and (\ref{eq:rateN2}):
\begin{align}
\frac{dE[N_{j}]}{dt}  &  =\frac{1}{n}\sum_{l=1}^{r}E\left[  (n_{j}-N_{j}%
)N_{l}\right]  Y_{w,l}^{T}Y_{h,j}-E[N_{j}]Y_{\delta,j}\nonumber\\
&  =\frac{1}{n}\sum_{l=1}^{r}(n_{j}-E[N_{j}])E[N_{l}]Y_{w,l}^{T}%
Y_{h,j}-E[N_{j}]Y_{\delta,j}-\frac{1}{n}\sum_{l=1}^{r}\mathrm{Cov}(N_{j}%
,N_{l})Y_{w,l}^{T}Y_{h,j}\nonumber\\
&  \approx\frac{1}{n}\sum_{l=1}^{r}(n_{j}-E[N_{j}])E[N_{l}]Y_{w,l}^{T}%
Y_{h,j}-E[N_{j}]Y_{\delta,j}-\frac{1}{n}\sum_{l=1}^{r}\Sigma_{\infty
,jl}Y_{w,l}^{T}Y_{h,j}. \label{corr_rate_eq}%
\end{align}
We find the metastable expectation by equating the time derivative to zero
($\frac{dE[N_{j}]}{dt}=0$). Hence, the corrected metastable state expectation
$E[N_{j}^{\infty}]=\widehat{N}_{j}^{\infty}$ in cluster $j$ is determined by
\[
\forall1\leq j\leq r:\frac{1}{n}\sum_{l=1}^{r}(n_{j}-\widehat{N}_{j}^{\infty
})\widehat{N}_{l}^{\infty}Y_{w,l}^{T}Y_{h,j}-\widehat{N}_{j}^{\infty}%
Y_{\delta,j}-\frac{1}{n}\sum_{l=1}^{r}\Sigma_{\infty,jl}Y_{w,l}^{T}Y_{h,j}=0.
\]
The expectation of a few nodes might turn out to be slightly negative with
this method, so we replace those by zero to obtain the corrected mean-field
approximation. It turns out that this correction is negligible (compared to
the original $N^{\infty}$) when using a low-dimensional $k$ factorisation
approximation of $\widetilde{A}$. However, in this example (using $r=n$ and
$k=n$), we found the following result:
\[%
\begin{tabular}
[c]{|c|c|c|c|}\hline
& Simulated & Mean-field & Corrected Mean-field\\\hline
Expectation & 1105.8 & 1136.1 & 1104.0\\\hline
Standard Deviation & 29.4 & 27.7 & 27.7\\\hline
\end{tabular}
\ \ \ \ \
\]
By construction, the mean field expectation is still about the same, but the corrected version is much more accurate. We also tried to iteratively update $\Sigma_\infty$ using the corrected expectation, but this did not give substantial improvements. Figure \ref{AirWeight_cdf} shows the good fit of the normal distribution, while Figure \ref{AirWeight_exp} assesses the quality of the corrected
mean-field approximation against simulations. We are not only able to substantially improve NIMFA, but we are
also able to obtain a rather precise indication of variances and covariances.
We ought to mention that our correction method is infeasible for very large networks;
there we would need to approximate $\widetilde{A}$ by a low rank non-negative
matrix factorisation. We advise to take $k$ as high as possible, and to take
the number of clusters $r$ as high as possible.

\section{Discussion}

Our analysis resembles a reaction-diffusion process studied in meta-populations
as outlined in \cite[Section IX]{PVM_RMP_epidemics2014}. Instead of
immediately assuming a mean-field approach as in reaction-diffusion models
\cite[Section IX]{PVM_RMP_epidemics2014}, our \textquotedblleft clustered
SIS\textquotedblright\ model makes the following approximations. First, we
assume that the infection rate matrix can be factorized as $\widetilde
{A}_{n\times n}=\left(  W^{T}\right)  _{n\times k}H_{k\times n}$ with a
tuneable $k$ (which can incur an approximation \cite{Hannah_LAA1983} if $k$ is
small compared to $n$). Second, we assume that $r$ clusters can be found for
which the cluster center represents the infectious state ($Z$-vector) for all
nodes in the same cluster well. After these two assumptions (that may lead to
approximations), a discrete Markov process $N$ arises, whose metastable state
can be obtained. Moreover, that metastable state $N^{\infty}$ exists provided
that the largest eigenvalue $\lambda$ of $\widetilde{A}^{T}/\delta$ exceeds 1,
precisely agreeing with earlier mean-field approximations (NIMFA) in
heterogeneous epidemics. Third, by additionally approximating the discrete
Markov chain $N$ by a continuous Brownian motion, a stochastic, non-linear
vector differential equation (\ref{SDE_non_linear_D}) is deduced, whose
linearization (which is a fourth approximation and presented in Appendix
\ref{sec_linearizing_SVDE}) for large population sizes can be solved exactly,
from which a Gaussian vector distribution (with mean and covariances) of the
number of infected nodes per cluster is obtained.

In spite of several approximations in the \textquotedblleft clustered
SIS\textquotedblright\ model, its accuracy, assessed in Section
\ref{sec_evaluation}, is remarkably good in synthetic graphs even for low dimensional factorizations.
These low dimensional factorizations are less sharp in a real-world airport network, but bypassing the factorization and many clusters, we were also able to predict the meta-stable distribution very well. Motivated to increase the overall
performance of the \textquotedblleft clustered SIS\textquotedblright\ model,
we found an efficient method to adjust the NIMFA bias by utilizing the
asymptotic covariance matrix $\Sigma_{\infty}$, defined in
(\ref{eq:Sigma_infty}) and computed in (\ref{Sigma_inf_explicit_Sigma}). In
essence, rather than neglecting pair-correlations as in NIMFA, we approximate
in SIS governing equation (\ref{governing_eq_heterogeneous_SIS}) the
pair-correlations by the asymptotic covariance matrix $\Sigma_{\infty}$. This new
correction method seems to be surprisingly accurate, which, we hope, will
inspire others to test and use it further in the future.

\medskip\textbf{Acknowledgements.} We are grateful to Huijuan Wang and Bo Qu
for the construction of the weighted airport network in Section
\ref{sec_application_real_world_networks} from the dataset of
openflights\footnote{https://openflights.org/data.html}.

{\footnotesize
\bibliographystyle{unsrt}
\bibliography{cac,MATH,misc,net,pvm,QTH,tel}
}

\newpage

\appendix{}

\section{Linearising the SVDE (\ref{SDE_non_linear_D}) for large $n$}

\label{sec_linearizing_SVDE}Each node is characterised by its infectiousness,
susceptibility and healing rate, captured by the vector $Z_{i}$ in
(\ref{def_Zi}). As the number of nodes $n$ increases, the average infection
rate between two nodes has to decrease, in order for a reasonable metastable
state to exist (otherwise all nodes will be infected at almost any time). An
asymptotic analysis for large $n$ requires us to describe how the network
grows with $n$. Therefore, we suppose that we can model the vectors $Z_{i}$ as
a sample of size $n$ from a distribution $\mu$ on ${\mathbb{R}}^{2k+1}$:
\[%
\begin{pmatrix}
\sqrt{n}W_{i}\\
\sqrt{n}H_{i}\\
\delta_{i}%
\end{pmatrix}
\sim\mu.
\]
Since our new approximate \textquotedblleft clustered SIS\textquotedblright%
\ model only looks at the inner products of $W_{i}$ and $H_{j}$, we can always
ensure that $W$ and $H$ \textquotedblleft live\textquotedblright\ at the same
scale. Now, we partition ${\mathbb{R}}_{+}^{2k+1}$ into $r$ disjoint sets
$\Upsilon_{1},\ldots,\Upsilon_{r}$ of positive $\mu$-mass, such that their
\textquotedblleft centers of mass\textquotedblright\ have small average
distance to the other points in the set. If we denote
\[
\rho_{j}=\mu(\Upsilon_{j})\mbox{ and }\tilde{Y}_{j}=\frac{1}{\rho_{j}}%
\int_{\Upsilon_{j}}{y}\mu(dy).
\]
then our clusters $C_{1},\ldots,C_{r}$ are defined by
\[
i\in C_{j}\Longleftrightarrow%
\begin{pmatrix}
\sqrt{n}W_{i}\\
\sqrt{n}H_{i}\\
\delta_{i}%
\end{pmatrix}
\in\Upsilon_{j}.
\]
Approximating $n_{j}/n\approx\rho_{j}$ yields
\[
Y_{w,j}=\frac{\sqrt{n}}{n_{j}}\sum_{i\in C_{j}}W_{i}\approx\tilde{Y}%
_{w,j}\mbox{  and  }Y_{h,j}=\frac{\sqrt{n}}{n_{j}}\sum_{i\in C_{j}}%
H_{i}\approx\tilde{Y}_{h,j}%
\]
and
\[
Y_{\delta,j}=\frac{1}{n_{j}}\sum_{i\in C_{j}}\delta_{i}\approx\tilde
{Y}_{\delta,j}.
\]
which shows that the matrix $Y_{h}^{T}Y_{w}$ tends to $\widetilde{Y}_{h}%
^{T}\widetilde{Y}_{w}$ when $n\rightarrow\infty$ by the strong law of large
numbers \cite[p. 119]{PVM_PAComplexNetsCUP}, since $n=\sum_{j=1}^{r}n_{j}$.

We return to the metastable state equation \eqref{eq:defNinfty} and define the
$r\times1$ vector $\tilde{N}^{\infty}\in\lbrack0,\rho_{1}]\times\ldots
\times\lbrack0,\rho_{r}]$ by the equations
\[
\forall\ j\in\{1,\ldots,r\}\ :\ \sum_{l=1}^{r}(\rho_{j}-\tilde{N}_{j}^{\infty
})\tilde{N}_{l}^{\infty}\tilde{Y}_{w,l}^{T}\tilde{Y}_{h,j}=\tilde{N}%
_{j}^{\infty}\tilde{Y}_{\delta,j}.
\]
Dividing the metastable state equation \eqref{eq:defNinfty} by $n$ shows that
for large $n$,
\[
N^{\infty}\approx n\tilde{N}^{\infty}.
\]
A closer look at the SVDE (\ref{SDE_non_linear_D}) reveals that we only need
to consider the matrix $M$ in (\ref{def_M}) as $n\rightarrow\infty$:
\begin{align*}
M  &  =\frac{1}{n}\mathrm{diag}(\mathbf{n}-N^{\infty})Y_{h}^{T}Y_{w}-\frac
{1}{n}\mathrm{diag}(Y_{h}^{T}Y_{w}N^{\infty})\\
&  \approx\mathrm{diag}(\mathbf{\rho}-\tilde{N}^{\infty})\tilde{Y}_{h}%
^{T}\tilde{Y}_{w}-\mathrm{diag}(\tilde{Y}_{h}^{T}\tilde{Y}_{w}\tilde
{N}^{\infty}),
\end{align*}
illustrating that $M$ also stabilises or converges. Thus, as $n\rightarrow
\infty$, the SVDE (\ref{SDE_non_linear_D}) reduces to
\begin{equation}
dD(t)=(M-\mathrm{diag}(Y_{\delta}))D(t)dt+\mathrm{diag}\left(  \sqrt
{2\mathrm{diag}(Y_{\delta})N^{\infty}/n}\right)  dB(t). \label{eq:linSDE}%
\end{equation}
This resulting SVDE (\ref{eq:linSDE}) is a standard linear SVDE in $D$, with a
constant diagonal covariance factor!

\subsection{Exact solution of SVDE (\ref{eq:linSDE})}

\label{sec_exact_solution_SVDE}Here, we demonstrate that the SVDE
(\ref{eq:linSDE}) can be solved exactly. After defining the $r\times r$
matrices 
\[
K=M-\mathrm{diag}(Y_{\delta})\mbox{ and }\Sigma=\mathrm{diag}\left(
2\mathrm{diag}(Y_{\delta})N^{\infty}\right)  /n,
\]
from (\ref{def_M}) and (\ref{eq:Ninfty}), the SVDE (\ref{eq:linSDE}) transforms into
\begin{equation}
dD(t)=KD(t)dt+\sqrt{\Sigma}\,dB(t), \label{eq:linSDE2}%
\end{equation}
and we obtain \cite{Harrison}
\[
D(t)=e^{tK}D(0)+\int_{0}^{t}e^{(t-s)K}\sqrt{\Sigma}\,dB(s).
\]
Since we are interested in the metastable solution, the eigenvalues of $K$
must have negative real part, so that
\[
\lim_{t\rightarrow\infty}e^{tK}D(0)=0.
\]
Also, the asymptotic covariance will be given by
\begin{equation}
\Sigma_{\infty}=\lim_{t\rightarrow\infty}\mathrm{Cov}\left[  D(t)\right]
=\int_{0}^{\infty}e^{sK}\Sigma\left(  e^{sK}\right)  ^{T}\,ds.
\label{eq:Sigma_infty}%
\end{equation}
There is another way to determine $\Sigma_{\infty}$: when we are in the
metastable state (and $n\rightarrow\infty$), the distribution of $D(t)$ should
not depend anymore on $t$. Thus, we may suppose that the metastable state is
given by
\[
D(t)\sim N(\mu,\Sigma_{\infty}).
\]
Invoking \eqref{eq:linSDE} shows that $D(t+dt)=D(t)+dD(t)$ also has a normal
distribution. Furthermore,
\[
E\left[  dD(t)\right]  =K\mu dt+\sqrt{\Sigma}\,E\left[
dB(t)\right]  =K\mu dt.
\]
So $\mu=0$ implies that the expectation is stable. Furthermore, since $dB(t)$
is independent of the past, we find
\begin{align*}
\mathrm{Cov}\left[  D(t+dt)\right]   &  =\mathrm{Cov}\left[
D(t)+KD(t)dt+\sqrt{\Sigma}\,dB(t)\right] \\
&  =\mathrm{Cov}\left[  (I+Kdt)D(t)\right]  +\Sigma dt\\
&  =(I+Kdt)\Sigma_{\infty}(I+K^{T}dt)+\Sigma dt\\
&  =\Sigma_{\infty}+(K\Sigma_{\infty}+\Sigma_{\infty}K^{T}+\Sigma
)dt+O(dt^{2}).
\end{align*}
In the metastable state, $\mathrm{Cov}\left[  D(t+dt)\right]  =\mathrm{Cov}%
\left[  D(t)\right]  $. Metastability thus implies that
\begin{equation}
K\Sigma_{\infty}+\Sigma_{\infty}K^{T}+\Sigma=0.
\label{metastability_matrix_equation}%
\end{equation}
One may verify\footnote{Indeed, since $Ke^{sK}=e^{sK}K$ commute (expand
$e^{sK}$ in a Taylor series), introducing $\Sigma_{\infty}$ in
(\ref{eq:Sigma_infty}) into (\ref{metastability_matrix_equation})
\begin{align*}
\int_{0}^{\infty}Ke^{sK}\Sigma\left(  e^{sK}\right)  ^{T}\,+e^{sK}%
\Sigma\left(  Ke^{sK}\right)  ^{T}ds+\Sigma &  =\int_{0}^{\infty}\frac{d}%
{ds}\left(  e^{sK}\Sigma\left(  e^{sK}\right)  ^{T}\right)  \,ds+\Sigma\\
&  =-\Sigma+\Sigma=0
\end{align*}
} that $\Sigma_{\infty}$, defined by \eqref{eq:Sigma_infty}, satisfies the
metastability matrix equation. To compute $\Sigma_{\infty}$, we will solve the
equation (\ref{metastability_matrix_equation}), rather than using
\eqref{eq:Sigma_infty}. Assuming that $K$ is diagonalizable, $K=V\Lambda
V^{-1}$ for some diagonal matrix $\Lambda$. Substituting $K=V\Lambda V^{-1}$
in (\ref{metastability_matrix_equation}), left-multiplying by $V^{-1}$ and
right-multiplying by $V^{-T}$ gives
\[
\Lambda V^{-1}\Sigma_{\infty}V^{-T}+V^{-1}\Sigma_{\infty}V^{-T}\Lambda
=-V^{-1}\Sigma V^{-T}.
\]
Since $\Lambda$ is diagonal, the $i,j$th coordinate satisfies
\[
(\Lambda_{ii}+\Lambda_{jj})\cdot(V^{-1}\Sigma_{\infty}V^{-T})_{ij}%
=(-V^{-1}\Sigma V^{-T})_{ij},
\]
and therefore, letting $J$ be the all-ones matrix,
\begin{equation}
\Sigma_{\infty}=-V\frac{V^{-1}\Sigma V^{-T}}{\Lambda J+J\Lambda}V^{T},
\label{Sigma_inf_explicit_Sigma}%
\end{equation}
where the division is elementwise.

We conclude that the metastable distribution $N=N^{\infty}+\sqrt{n}D$ will be
an $r$-dimensional normal distribution, with expectation equal to $N^{\infty}$
(since $E\left[  D\right]  =0$), and covariance matrix
$n\Sigma_{\infty}$.

\end{document}